\documentclass[12pt]{article}
\usepackage{amsmath}
\usepackage{amssymb}
\usepackage{graphicx}
\usepackage[cp1251]{inputenc}
\usepackage[russian]{babel}
\newcommand{\il}[2]{\int\limits_{#1}^{#2}}

\newcommand{\ph}{\phantom{a}}
\newcommand{\phh}{\phantom{aaa}}

\newcommand{\sist}[2]{\left\{
\begin{array}{l}
{#1}\\
\ph\\
{#2}
\end{array}
\right.}

\begin{document}

\vskip 20pt

MSC 34C10

\vskip 20pt

\centerline{\bf A non oscillation criterion for}
 \centerline{\bf extended quasi linear Hamiltonian systems}

\vskip 20 pt

\centerline{\bf G. A. Grigorian}
\centerline{\it Institute  of Mathematics NAS of Armenia}
\centerline{\it E -mail: mathphys2@instmath.sci.am}
\vskip 20 pt

\noindent
Abstract. A new approach with the Riccati equation method is used to obtain a  non oscillation criterion for extended quasi linear Hamiltonian systems.

\vskip 20 pt


\vskip 20 pt

{\bf 1. Introduction}.
Let $A(t,u,v), \ph B(t,u,v)$ and $C(t,u,v)$ be complex-valued locally integrable in $t$  and continuous in $u$ and $v$ matrix  functions of dimension $n \times n$ on $[t_0,+\infty)\times\mathbb{C}^n\times\mathbb{C}^n$ and let $B(t,u,v)$ be Hermitian,  i.e.  $B(t,u,v) = P^*(t,u,v), \ph t\ge t_0, \ph u,v \in \mathbb{C}^n$, where $*$ denotes the conjugation sign. Let $\mu(t,u,v)$ be real-valued locally integrable in $t$ continuous in $u$ and $v$ scalar function on $[t_0,+\infty)\times\mathbb{C}^n\times\mathbb{C}^n$.
Consider the generalized non linear Hamiltonian system
$$
\sist{\phi' = A(t,\phi,\psi)\phi + B(t,\phi,\psi) \psi,}{\psi' = C(t,\phi,\psi)\phi + [\mu(t,\phi,\psi) I - A^*(t,\phi,\psi)]\psi, \ph t\ge t_0,} \eqno (1.1)
$$
where $I$ is the identity matrix of dimension $n\times n$. The non linear (in particular linear) Hamiltonian systems play important role in  many areas of natural sciences: in  physics and mechanics, in particular, in many-dimensional oscillator systems, in the particle dynamics,  chaos and phase space diffusion,  quantum mechanics, in applied mathematics: in variational theory, optimal control and filtering, dynamic programming and differential games, invariant embedding  and scattering processes (see [2] and cited works therein). To the study of nonlinear as well as linea Hamiltonian systems are devoted many works (see e.g. [1,3--8] and cited works therein).  By a solution of the system (1.1) on any interval $[t_1,t_2) \subset [t_0,+\infty)$ we mean an ordered pair $(\phi(t),\psi(t))$ of absolutely continuous vector functions $\phi(t)$ and $\psi(t)$  of dimension $n$ on $[t_1,t_2)$, satisfying (1.1) almost everywhere on $[t_1,t_2)$. According to general theory of normal systems of ordinary differential equations for every $t_1 \ge t_0$ and for every initial values $\phi_0, \ph \psi_0\in\mathbb{C}^n$ there exists $t_2 > t_1$ such that the system  (1.1) has a solution $(\phi(t),\psi(t))$ on $[t_1,t_2)$ with $\phi(t_1) = \phi_0, \ph \psi(t_1) = \psi_0$. Main interest from the point of view of qualitative theory of differential equation represents the case $t_2=+\infty$. In this case one says  that the system (1.1) is global solvable or has a global solution. Conditions, providing the global solvability of the system (1.1), are the following quasi linearity conditions
$$
\alpha) \ph ||A(t,u,v)|| \le a_0(t), \ph ||B(t,u,v)|| \le b_0(t), \ph ||C(t,u,v)|| \le c_0(t), \ph |\mu(t,u,v)| \le \mu_0(t),
$$
$t \ge t_0,$
where for any square matrix $M$ the symbol $||M||$ denotes any euclidian norm of $M$, $a_0(t), \ph b_0(t), \ph c_0(t), \ph \mu_0(t)$ are locally integrable functions on $[t_0,+\infty)$. Under these conditions the Cauchy problem
$$
\phi(t_1) = \phi_1, \ph \psi(t_1) =\psi_1, \ph t_1 \ge t_0, \ph \phi_1,\psi_1 \in \mathbb{C}^n.
$$
for the system (1.1)  has the unique solution on $[t_0,+\infty)$, which can be proved using the contracting mapping principle (see [5]).

\vskip 10pt

{\bf Definition 1.1.} {\it A solution $(\phi(t),\psi(t))$ of the system (1.1), existing  on $[t_0,+\infty)$, is called conjoined, if there exists $t_1\ge t_0$ such that $\phi(t_1) \ne \theta\equiv(0,\ldots,0)$ and $\psi(t_1) = Y\phi(t_1),$ where $Y$ is a matrix of dimension $n\times n$ with the nonnegative definite  Hermitian part (i.~ e.  $Y + Y^*$ is nonnegative definite).
}

\vskip 10pt

{\bf Definition 1.2.} {\it A solution $(\phi(t),\psi(t))$ of the system (1.1), existing on $[t_0,+\infty)$, is called non oscillatory if there exists $T\ge t_0$ such that $\phi(t) \ne \theta, \ph t \ge T.$
}

\vskip 10pt

{\bf Definition 1.3.} {\it The system (1.1) is called nonsingular if it has a conjoined solution on $[t_0,+\infty)$.}

\vskip 10pt

{\bf Definition 1.4.} {\it The system (1.1) is called non oscillatory if it is nonsingular and every its conjoined solution is non oscillatory.
}

 Note that the system (1.1) with the restrictions $\alpha)$ is nonsingular. The question of existence of nonsingular systems, different from quasi linear ones, remains open.

In this paper we use the Riccati equation method to prove an non oscillation criterion for the system (1.1).

\vskip 10pt

{\bf 2. Main result.} The nonnegative (positive) definiteness of any square matrix $M$ is usually denoted by $M\ge 0 \ph (>0)$.

\vskip 10pt

{\bf Theorem 2.1.} {\it Let the conditions $\alpha)$ be satisfied. If $$B(t,u,v) \ge 0, \ph C(t,u,v) + C^*(t,u,v) \ge 0, \ph t \ge t_0,\ph u,v\in\mathbb{C}^n,$$ then the system (1.1) is non oscillatory.
}
\phantom{aaaaaaaaaaaaaaaaaaaaaaaaaaaaaaaaaaaaaa}$\blacksquare$

{\bf Remark 2.1.} {\it It will be clear from the proof of Theorem 2.1 that the conditions $\alpha)$ can be replaced by others (if they exist), providing the non singularity of the system (1.1).

}

\vskip 5pt

{\bf 3. Auxiliary propositions}. The next three lemmas have importance  in the proof of the main result. Let $M_l \equiv (m_{ij}^l)_{i,j=1}^n, \ph l=1,2$ be complex-valued matrices.

{\bf Lemma 3.1.} {\it The equality
$$
tr(M_1 M_2) = tr (M_2 M_1)
$$
is valid.}

Proof.  We have $tr (M_1 M_2) = \sum\limits_{j=1}^n(\sum\limits_{k=1}^n m_{jk}^1 m_{kj}^2) = \sum\limits_{k=1}^n(\sum\limits_{j=1}^n m_{jk}^1 m_{kj}^2) = \sum\limits_{k=1}^n(\sum\limits_{j=1}^n m_{kj}^2 m_{jk}^1) = tr (M_2 M_1).$ The lemma is proved.

{\bf Lemma 3.2.} {\it Let $H_j, \ph j=1,2$ be Hermitian matrices such that $H_j\ge 0, \ph j=1,2.$ Then
$$
tr (H_1 H_2) \ge 0.
$$
}

Proof. Let $U$ be an unitary transformation such that $\widetilde{H}_1 \equiv U H_1 U^* = diag \{h_1,\dots,h_n\}$.
Since any unitary transformation preserves the nonnegative definiteness of any Hermitian  matrix    we have
$$
h_j \ge 0, \phh j=\overline{1,n}. \eqno (3.1)
$$
Let $\widetilde{H}_2 \equiv U H_2 U^* = (h_{ij})_{i,j=1}^n.$ As for as $\widetilde{H}_2$ is Hermitian it follows that (see [4], p. 300, Theorem 20)  $h_{jj} \ge 0, \ph j=\overline{1,n}$. This together with (3.1) implies
$$
tr (H_1 H_2) = tr([U H_1 U^*] [U H_2 U^*]) = tr (\widetilde{H}_1 \widetilde{H}_2) = \sum\limits_{j=1}^n h_j h_{jj} \ge 0.
$$
The Lemma is proved.

{\bf Lemma 3.3.} {\it Let $H \ge 0$ be a Hermitian matrix of dimension $n \times n$ and let $V$ be any matrix of the same dimension. Then
$$
V H V^* \ge 0. \eqno (3.4)
$$
}

Proof. For any vectors $x$ and $y$ of dimension $n$ denote by $\langle x,y\rangle$ their scalar product. Then $\langle V H V^* x, x\rangle = \langle H(V^* x), (V^* x)\rangle \ge 0$ (since $H \ge 0$). Hence (3.4) is valid. The lemma is proved.

\vskip 10pt

{\bf 4. Proof of Theorem 2.1.} As was noted in the introduction, if the conditions $\alpha)$ are satisfied, the system (1.1) becomes non-singular. Let then $(\phi_0(t),\psi_0(t))$ be a conjoined solution of the system (1.1), and let $\psi(t_1) =Y_0\phi(t_1)$ for some $t_1 \ge t_0$,  where $Y_0$ is a matrix of dimension $n\times n$, such that $Y_0 + Y_0^*\ge 0$.
Consider the matrix Riccati equation
$$
Y' + Y B_0(t) Y + Y A_0(t) + [A_0^*(t)  -\nu_0(t) I]Y -C_0(t) = 0, \ph t \ge t_1, \eqno (4.1)
$$
where $A_0(t)\equiv A(t,\phi_0(t),\psi_0(t)), \ph B_0(t)\equiv B(t,\phi_0(t),\psi_0(t)), \ph C_0(t)\equiv C(t,\phi_0(t),\psi_0(t)), \linebreak \nu_0(t)\equiv \mu(t,\phi_0(t),\psi_0(t)), \ph t \ge t_1.$ Due to the conditions $\alpha)$ we have $||A_0(t)|| \le a_0(t), \linebreak ||B_0(t)|| \le b_0(t), \ph ||C_0(t)|| \le c_0(t), \ph |\nu_0(t)| \le \mu_0(t), \ph t \ge t_0$. Therefore the matrix functions $A_0(t), \ph B_0(t), \ph  C_0(t)$ and the scalar function $\nu_0(t)$ are locally integrable over $[t_0,+\infty)$. Then it is reasonable to mean a solution of Eq. (4.1) on any interval $[\tau_1,\tau_2) \subset [t_0,+\infty)$ an absolutely continuous on $[\tau_1,\tau_2)$ matrix function of dimension $n\times n$, satisfying (4.1) almost everywhere on $[\tau_1,\tau_2)$.  Let $Y(t)$ be the solution of Eq. (4.1) with $Y(t_1) = Y_0$, and let $[t_1,t_2)$ be its maximum existence interval. It is not difficult to verify that $Y(t)$ and $(\phi_0(t),\psi_0(t))$ are connected by equalities
$$
\phi_0'(t) = [A_0(t) + B_0(t) Y(t)]\phi_0(t), \ph \psi_0(t) = Y(t)\phi_0(t), \ph t \in [t_1,t_2).
$$
The first of the obtained equalities shows that we can interpret $\phi_0(t)$ as a nontrivial solution of the linear system
$$
\phi' = [A_0(t) + B_0(t) Y(t)]\phi, \ph t \in [t_1.t_2).
$$
Therefore, by the uniqueness theorem  the proof of the theorem will be completed, if we show that
$$
t_2 = +\infty. \eqno (4.2).
$$
By (4.1) we have
$$
Y'(t) + Y(t)B_0(t)Y(t) + Y(t) A_0(t) +[A_0^*(t) - \nu_0(t) I]Y(t) - C_0(t) = 0, \ph t \in [t_1,t_2),
$$
$$
[Y^*(t)]' + Y^*(t)B_0(t)Y^*(t) + A^*_0(t)Y^*(t) + Y^*(t) [A_0(t) - \nu_0(t) I] - C_0^*(t) = 0, \ph t \in [t_1,t_2),
$$
Summing up these equalities and making some simplifications we obtain
$$
(Y(t) + Y^*(t))' + (Y(t) + Y^*(t)) B_0(t)(Y(t) + Y^*(t)) + \Bigl[A_0^*(t) - \frac{\nu_0(t)}{2} I\Bigr](Y(t) + Y^*(t)) +
$$
$$
+(Y(t) + Y^*(t))\Bigl[A_0(t)- \frac{\nu_0(t)}{2} I\Bigr] - C_0(t) - C_0^*(t)  - Y(t)B_0(t)Y^*(t) - Y^*(t)B_0(t)Y(t) = 0, \eqno (4.3)
$$
$t\in[t_1,t_2).$ Assume $Y(t_1) + Y^*(t_1) > 0$. Show that
$$
Y(t) + Y^*(t) > 0, \ph t \in [t_1,t_2). \eqno (4.4)
$$
Suppose that this is not true. Then there exists $t_3\in (t_1,t_2)$ such that
$$
Y(t) + Y^*(t) > 0, \ph t \in [t_1,t_3) \eqno (4.5)
$$
and
$$
\det [Y(t_3) + Y^*(t_3)] = 0. \eqno (4.6)
$$
It follows from (4.5) that the matrix function $Y(t) + Y^*(t)$ does not degenerate on $[t_1,t_3)$, i.e., $(Y(t)+Y^*(t))^{-1}$ exists on $[t_1,t_3).$ Then by (4.3) one can interpret the matrix function $U(t) \equiv Y(t) + Y^*(t), \ph t \in [t_1,t_3)$ as a solution of the linear equation
$$
U' + \Bigl\{[Y(t) + Y^*(t)]B_0(t) + A_0^*(t) - \nu_0(t) I + [Y(t) + Y^*(t)]A_0(t)[Y(t) + Y^*(t)]^{-1} -
$$
$$
-[C_0(t) + C_0^*(t) + Y(t)B_0(t)Y^*(t) + Y^*(t)B_0(t)Y(t)](Y(t) + Y^*(t))^{-1}\Bigr\} U = 0, \ph t \in [t_1,t_3).
$$
Then by virtue of the Liouville  formula we have
$$
\det[Y(t) + Y^*(t)] = \det[Y(t_1) + Y^*(t_1)]\exp\biggl\{-\il{t_1}{t}\Bigl\{[Y(\tau) + Y^*(\tau)]B_0(\tau) +
 $$
 $$
 +A_0^*(\tau) - \nu_0(\tau) I + [Y(\tau) + Y^*(\tau)]A_0(\tau)[Y(\tau) + Y^*(\tau)]^{-1} -
$$
$$
-[C_0(\tau) + C_0^*(\tau) + Y(\tau)B_0(\tau)Y^*(\tau) + Y^*(\tau)B_0(\tau)Y(\tau)](Y(\tau) + Y^*(\tau))^{-1}\Bigr\}d\tau\biggr\},  \eqno (4.7)
$$
$t\in [t_1,t_3)$. By Lemma 3.3 $Y(t)B_0(t)Y^*(t) + Y^*(t)B_0(t)Y(t)\ge 0, \ph t\in [t_1,t_0)$. Then, since $(Y(t) + Y^*(t))^{-1}> 0,\ph t \in [t_1,t_3)$ by Lemma 3.2 $tr \Bigl[(Y(t)B_0(t)Y^*(t) + Y^*(t)B_0(t)Y(t))(Y(t)+Y^*(t))\Bigr] \ge 0, \ph t \in [t_1,t_3)$. Moreover, by Lemma 3.2 it follows from the condition $C(t,u,v) + C^*(t,u,v) \ge 0, \ph t \ge t_0, u,v \in\mathbb{C}^n$ that $tr\Bigl[(C_0(t) + C_0^*(t))(Y(t)+Y^*(t))^{-1}\Bigr] \ge 0, \ph t \in [t_1,t_3)$. Therefore, it follows from  (4.7) that
$$
|\det[Y(t) + Y^*(t)]| \ge|\det[Y(t_1) + Y^*(t_1)]|\ge e^{-c}|\det[Y(t_1) + Y^*(t_1)], \ph t\in[t_1,t_3],
$$
where $c\equiv\max\limits_{t\in[t_1,t_3]}\Bigl|\il{t_1}{t}tr\Bigl[(Y(\tau) + Y^*(\tau))B_0(\tau) +A_0(\tau) +A_0^*(\tau) -\nu_0(\tau)I\Bigr]d\tau\Bigr|$, which contradicts (4.6). The obtained contradiction proves (4.4).
Consider the linear system
$$
\sist{\Phi' = A_0(t)\Phi + B_0(t)\Psi,}{\Psi' = C_0(t)\Phi +[\nu_0(t)I - A_0^*(t)]\Psi, \ph t \ge t_1.} \eqno (4.8)
$$
By a solution of this system we mean an ordered pair $(\Phi(t),\Psi(t))$ of absolutely continuous matrix functions of dimension $n\times n$ on $[t_1,+\infty)$, satisfying (4.8) almost everywhere on $[t_1,+\infty)$. Substituting
$$
\Psi = Y \Phi, \ph t \ge t_1
$$
in (4.8) one can verify that all solutions $Y_1(t)$ of Eq. (4.1), existing on any interval $[\tau_1,\tau_2) \subset [t_1,+\infty)$, are connected with solutions $(\Phi(t),\Psi(t))$ of the system (4.8) by relations
$$
\Phi'(t) = [A_0(t) + B_0(t)Y_1(t)]\Phi(t), \phh \Psi(t) = Y_1(t)\Phi(t), \phh t \in [\tau_1,\tau_2). \eqno (4.9)
$$
By the Liouville formula from here we get
$$
\det \Phi(t) = \det \Phi(\tau_1)\exp\biggl\{\il{\tau_1}{t} tr\Bigl[A_0(\tau) + B_0(\tau)Y_1(\tau)]d\tau\biggr\}, \ph t \in [\tau_1,\tau_2),
$$
$$
\overline{\det \Phi(t)} = \overline{\det \Phi(\tau_1)}\exp\biggl\{\il{\tau_1}{t} tr\Bigl[A_0^*(\tau) + Y_1^*(\tau)B_0(\tau)]d\tau\biggr\}, \ph t \in [\tau_1,\tau_2).
$$
Hence, (since by  Lemma 3.1 $tr Y^*(t)B_0(t) = tr B_0(t)Y^*(t)$)
$$
|\det \Phi(t)|^2 = |\det \Phi(\tau_1)|^2\exp\biggl\{\il{\tau_1}{t} tr\Bigl[A_0(\tau) + A_0^*(\tau) + B_0(\tau)(Y_1(\tau)+Y_1^*(\tau))]d\tau\biggr\}, \eqno (4.10)
$$
$t\in [\tau_1,\tau_2)$. Let $(\Phi_1(t),\Psi_1(t))$ be a solution of the system (4.8) with $\Phi_1(t_1) = I, \ph \Psi_1(t_1) = Y(t_1).$ Then by (4.10)
$$
|\det \Phi_1(t)|^2 = |\det \Phi_1(t_1)|^2\exp\biggl\{\il{t_1}{t} tr\Bigl[A_0(\tau) + A_0^*(\tau) + B_0(\tau)(Y(\tau)+Y^*(\tau))]d\tau\biggr\}, \eqno (4.11)
$$
$t\in [t_1,t_2)$. According to a condition of the theorem $B(t,u,v) \ge 0, \ph t \ge t_0, \ph u,v \in \mathbb{C}^n$.
Hence, $B_0(t)\ge 0, \ph t \ge t_1$. By Lemma 3.2 it follows from here, from (4.4) and from (4.11) that  that
$$
|\det \Phi_1(t_2)|^2 \ge |\det \Phi_1(t_1)|^2\exp\biggl\{\il{t_1}{t_2} tr\Bigl[A_0(\tau) + A_0^*(\tau)]d\tau\biggr\} > 0.
$$
Hence, $\det \Phi_1(t) \ne 0, \ph t \in [t_1,t_2 + \varepsilon)$ for some $\varepsilon >0$. In virtue of (4.9) it follows from here that $Y_1(t)\equiv \Psi_1(t)\Phi_1^{-1}(t), \ph t \in [t_1,t_2+\varepsilon)$ is a solution of Eq. (4.1) on $[t_1,t_2+\varepsilon)$, and coincides with $Y(t)$ on $[t_1,t_2)$. By the uniqueness theorem it follows from here that $[t_1,t_2)$ is not the maximum existence interval for $y(t)$, which contradicts our supposition. The obtained contradiction proves (4.2). Thus we prove (4.2) for the particular case $Y(t_1) +~ Y^*(t_1) > 0$. Let us prove (4.2) in the general case  $Y(t_1) + Y^*(t_1) \ge 0$. For every $\delta > 0$ denote by $Y_\delta(t)$ the solution of Eq. (4.1), satisfying the initial condition $Y_\delta(t_1) = Y(t_1) + \delta I$. Obviously, $Y_\delta(t_1)+ Y_\delta^*(t_1) \ge 2\delta I > 0$. Then similar to the already proven one can show that $Y_\delta(t)$ exists on $[t_1,+\infty)$ and
$$
Y_\delta(t) + Y_\delta^*(t) > 0, \phh t \ge t_1. \eqno (4.12)
$$
Let $\lambda(t)$ be the least eigenvalue for $Y(t) + Y^*(t)$ on $[t_1,t_2)$ and let $\lambda_\delta(t)$ be the least eigenvalue of $Y_\delta(t) + Y_\delta^*(t)$ on $[t_1,+\infty)$. It follows from (4.12) that
$$
\lambda_\delta(t) > 0, \phh t \ge t_1. \eqno (4.13)
$$
Since the solutions of Eq. (4.1) are continuously dependent on their initial values we have $\lim\limits_{\delta\to 0+}\lambda_\delta(t) = \lambda(t), \ph t \ge t_1.$ This together with (4.13) implies that  $\lambda(t) \ge 0, \ph t\in[t_1,t_2)$. Therefore, $Y(t)+Y^*(t) \ge 0, \ph t \in [t_1,t_2)$. Further the proof of (4.2) is as the proof of (4.2)
in the particular case $Y(t_1) + Y^*(t_1) > 0$. The proof of the theorem is completed.

\vskip 20 pt

\centerline{ \bf References}

\vskip 20pt

\noindent
1. K. I. Al-Dosary, H. K. Abdullah and D. Hussein, Short Note on  oscillation of matrix \linebreak \phantom{a} Hamiltonian systems. Yokohama Math. J. vol 50, 2003, pp. 23--30.

\noindent
2. P. M. Cincotta, Notes on
Non-Linear Hamiltonian Dynamics (in
progress)
June 13, \linebreak \phantom{a} 2024, 281 pages.

\noindent
3. Sh. Chen, Z. Zheng, Oscillation criteria of Yan type for linear Hamiltonian systems, \linebreak \phantom{a} Comput.  Math. with Appli., 46 (2003), 855--862.

\noindent
4.  F. R. Gantmacher, Theory of Matrix. Second Edition (in Russian). Moskow,,\linebreak \phantom{a} ''Nauka'', 1966.

\noindent
5. G. A. Grigorian,
The Cauchy problem for quasilinear systems of functional  \linebreak \phantom{a} differential equations.
Sarajevo J. Math. 18(31) (2022), no. 2, 265–272.

\noindent
6. G. A. Grigorian, Oscillation criteria for linear matrix Hamiltonian systems. \linebreak \phantom{aa} Proc. Amer. Math. Sci, Vol. 148, Num. 8 ,2020, pp. 3407 - 3415.

\noindent
7. I. S. Kumary and S. Umamaheswaram, Oscillation criteria for linear matrix \linebreak \phantom{a} Hamiltonian systems, J. Differential Equ., 165, 174--198 (2000).

\noindent
8. L. Li, F. Meng and Z. Zheng, Oscillation results related to integral averaging technique\linebreak \phantom{a} for linear Hamiltonian systems, Dynamic Systems  Appli. 18 (2009), \ph \linebreak \phantom{a} pp. 725--736.

\end{document}